\input amstex
\documentstyle{amsppt}
\magnification=\magstep1
 \hsize 13cm \vsize 18.35cm \pageno=1
\loadbold \loadmsam
    \loadmsbm
    \UseAMSsymbols
\topmatter
\NoRunningHeads
\title Note on $q$-Dedekind type sums related to $q$-Euler
polynomials
\endtitle
\author
 Taekyun Kim
\endauthor

\abstract Recently, $q$-Dedekind type sums related to $q$-zeta
function and  basic $L$-series are studied by Simsek in [Y. Simsek,
$q$-Dedekind type sums related to $q$-zeta function and basic
$L$-series, J. Math. Anal. Appl. 318(2006) 333-351] and Dedekind
type sums related to Euler numbers and polynomials are introduced in
previous paper[T. Kim, Note on Dedekind type DC sums, Adv. Stud.
Contemp. Math. 18(2009), 249-260]. It is the purpose of this paper
to construct a $p$-adic continuous function for an odd prime to
contain a $p$-adic $q$-analogue of higher order Dedekind type sums
related to $q$-Euler polynomials and numbers by using an invariant
$p$-adic $q$-integrals.
\endabstract
\thanks  2000 AMS Subject Classification: 11B68, 11S80
\newline  The present Research has been conducted by the research
Grant of Kwangwoon University in 2008
\endthanks
\endtopmatter

\document

{\bf\centerline {\S 1. Introduction/ Preliminaries}}

 \vskip 15pt
Let $p$ be a fixed odd prime number. Throughout this paper $\Bbb
Z_p$, $\Bbb Q_p$, $\Bbb C$ and $\Bbb C_p$ will, respectively, denote
the ring of $p$-adic rational integers, the field of $p$-adic
rational numbers, the complex number field and the completion of
algebraic closure of $\Bbb Q_p$. Let $v_p$ be the normalized
exponential valuation of $\Bbb C_p$ with
$|p|_p=p^{-v_p(p)}=\frac{1}{p}.$ When one talks of $q$-extension,
$q$ is variously considered  as and indeterminate, a complex number
$q\in \Bbb C$ or $p$-adic number $q\in \Bbb C_p$. If $q\in \Bbb C$,
one normally assumes $|q|<1$. If $q\in\Bbb C_p$, one normally
assumes $|1-q|_p<1.$ Recently, we proposed a definition of a
$q$-extension of $p$-adic Haar measure as follows: For any positive
integer $N$, we set
$$\mu_q(a +p^N\Bbb Z_p)=\frac{(1+q)(-q)^a}{1+q^{p^N}}, \text{ (see [1-12]).} $$
for $0\leq a \leq p^N -1$ and this can be extended to a measure on
$\Bbb Z_p$. This measure yields an invariant $p$-adic $q$-integral
for each non-negative integer $m$ and the $m$-th Carlitz's type
$q$-Euler numbers $ \varepsilon_{n, q}$ can be represented by this
$p$-adic $q$-integral as follows:
$$\varepsilon_{n, q}=\int_{\Bbb Z_p}\left(\frac{1-q^a}{1-q}\right)^m
d\mu_q(a)=\lim_{N\rightarrow
\infty}\frac{1+q}{1+q^{p^N}}\sum_{x=0}^{p^N-1}\left(\frac{1-q^a}{1-q}\right)^m(-q)^a,$$
which has  a sense  as we see readily that the limit is convergent
(see[8,9]). The modified $q$-Euler numbers are also defined as
$$E_{m,q}=\int_{\Bbb Z_p}\left(\frac{1-q^a}{1-q}\right)^m q^{-a}
d\mu_q(a), \text{ (see [6] )}.$$ Note that $ \lim_{ q\rightarrow
1}E_{n,q}=E_n ,$ where $E_n$ are the $n$-th Euler numbers. Now, we
also consider  the $q$-Euler polynomials as follows:
$$E_{m,q}(x)=\int_{\Bbb Z_p}q^{-t}
\left(\frac{1-q^{x+t}}{1-q}\right)^m d\mu_q(t), \text{ for $x\in\Bbb
Z_p, m\in\Bbb N.$}$$ These numbers $E_{m,q}(x)$ can be represented
by
$$E_{m,q}(x)=\sum_{l=0}^m
\binom{m}{l}q^{xl}E_{l,q}\left(\frac{1-q^x}{1-q}\right)^{m-l}. $$
For any positive integer $h, k$ and $m$, Dedekind type DC sums are
defined as
$$ S_m( h, k)=\sum_{M=1}^{k-1}(-1)^{M-1}\frac{M}{k} \bar E_m ( \frac{hM}{k}), \text{ (see [7]),} $$
where $ \bar E_m(x)$ are the $m$-th periodic Euler function.

By using an invariant $p$-adic $q$-integral on $\Bbb Z_p$, we
construct a $p$-adic continuous function for an odd prime to contain
a $p$-adic $q$-analogue of higher order Dedekind type DC sums $k^m
S_{m+1}(h,k)$ in this paper.  It is the purpose of this paper  to
give a $q$-analogue of $p$-adic Dedekind type DC sums by using
invariant $p$-adic $q$-integral on $\Bbb Z_p$ approach the $p$-adic
analogue of higher order Dedekind sums  at $q=1$ as follows:

\proclaim{ Theorem} Let $h, k$ be positive integer with $(h, k)=1$,
$p\nmid k$. For $s \in \Bbb Z_p$, let us define $p$-adic Dedekind
type DC sums as follows.
$$
S_{p,q}(s: h, k: q^k)=
\sum_{M=1}^{k-1}\left(\frac{1-q^M}{1-q}\right)(-1)^{M-1}T_q(s, hM,
k: q^k).$$ Then there exists a continuous function $S_{p,q}(s:h,
k:q^k) $ on $\Bbb Z_p$ which satisfies
$$\aligned
&S_{p,q}(m:h, k:q^k)=\left(\frac{1-q^k}{1-q}\right)^{m+1}S_{m,q}(h,
k:q^k)\\
&-\left(\frac{1-q^k}{1-q}\right)^{m+1}\left(\frac{1-q^{kp}}{1-q^k}\right)^{m}S_{m,q}(
(p^{-1}h)_k, k: q^{pk}),
\endaligned$$
 where $ m+1\equiv 0$ (mod $p-1$), and $(p^{-1}a)_N$ denotes  the
 integer $x$ with $0\leq x <N$, $ px\equiv a$ ( mod $N$ ).

 \endproclaim

\vskip 10pt

{\bf\centerline {\S 2. Proof of theorem}} \vskip 10pt

The $q$-Euler numbers $E_{m, q}$ can be written as
$$ E_{0,q}=\frac{1+q}{2},\qquad ( qE+1)^n +E_{n,q}=0 \text{ if $n>0$},$$
which is
$$E_{n,q}=(1+q)\left( \frac{1}{1-q}\right)^n \sum_{l=0}^n
\binom{n}{l} (-1)^l \frac{1}{1+q^l}, \text{ (see [6]),} $$ where we
use the technique method notation by replacing $E^n$ by $E_{n,q},$
symbolically. Let $w$ denote Teich$\ddot{u}$ller character (mod
$p$). For $ x\in \Bbb Z_p^{\times}=\Bbb Z_p\backslash p\Bbb Z_p$, we
set $$<x>_q=<x:q>=w^{-1}(x)\left(\frac{1-q^x}{1-q}\right).$$ Let $a$
and $N$ be positive integers with $(p,a)=1$ and $p|N. $ Define
$$T_q(s, a, N:
q^N)=w^{-1}(a)<a>_q^s
\sum_{j=0}^{\infty}\binom{s}{j}q^{aj}\left(\frac{1-q^N}{1-q^a}\right)^jE_{j,
q^N}, \text{ for $s\in\Bbb Z_p$}.$$ In particular , if $ m+1\equiv 0
 $ (mod $p-1$), then
$$\aligned
T_q(m,a, N:
q^N)&=\left(\frac{1-q^a}{1-q}\right)^m\sum_{j=0}^m\binom{m}{j}q^{aj}\left(\frac{1-q^N}{1-q^a}\right)^jE_{j, q^N} \\
&=\left(\frac{1-q^N}{1-q}\right)^m\int_{\Bbb Z_p}\left(
\frac{1-q^{Nx+a}}{1-q^N}\right)^m q^{-Nx} d\mu_{q^N}(x).
\endaligned$$
 Therefore, $T_q(m, a, N: q^N)$ is  a continuous $p$-adic  extension
 of $\left(\frac{1-q^N}{1-q}\right)^mE_{m, q^N}(\frac{a}{N}).$

Let $[ \cdot ]$ be Gauss' symbol and let $\{ x \}=x-[x].$ Then we
now consider a $q$-analogue of  higher order Dedekind type DC sums
$S_{m, q}(h, k:q^l)$ as follows:
$$S_{m,q}(h,k:q^l)
=\sum_{M=1}^{k-1}(-1)^{M-1}\left(\frac{1-q^M}{1-q^k}\right)\int_{\Bbb
Z_p}q^{-lx}\left(\frac{1-q^{l(x+\{\frac{hM}{k}\})}}{1-q^l}\right)^m
d\mu_{q^l}(x).$$ If $m+1\equiv 0$ (mod$p-1$), then we have
$$\aligned
&\left(
\frac{1-q^k}{1-q}\right)^{m+1}\sum_{M=1}^{k-1}\left(\frac{1-q^M}{1-q^k}\right)(-1)^{M-1}\int_{\Bbb
Z_p}\left(\frac{1-q^{k(x+\frac{hM}{k})}}{1-q^k} \right)^m
q^{-kx}d\mu_{q^{-k}}(x)\\
&=\sum_{M=1}^{k-1}(-1)^{M-1}\left(\frac{1-q^M}{1-q}\right)\left(\frac{1-q^k}{1-q}\right)^m\int_{\Bbb
Z_p} \left( \frac{1-q^{k(x+\frac{hM}{k})}}{1-q^k}\right)^mq^{-kx}
d\mu_{q^k}(x),
\endaligned\tag1$$
where $p|k, \,(hM, p)=1$ for each $M$. From (1), we note that
$$\aligned
&\left(\frac{1-q^k}{1-q}\right)^{m+1}S_{m,q}(h,k:q^k)\\
&=\sum_{M=1}^{k-1}\left(\frac{1-q^M}{1-q}\right)
\left(\frac{1-q^k}{1-q}\right)^m(-1)^{M-1}\int_{\Bbb
Z_p}q^{-kx}\left(\frac{1-q^{k(x+\frac{hM}{k})}}{1-q^k}\right)^m
d\mu_{q^k}(x)\\
&=\sum_{M=1}^{k-1}(-1)^{M-1}\left(\frac{1-q^M}{1-q}\right)T_q(m,
(hM)_k:q^k),
\endaligned\tag2$$
where $(y)_k$ denotes the integers $x$ such that $0\leq x <n$ and $
x\equiv \alpha $ ( mod $k$).

It is easy to check that
$$\aligned
&\int_{\Bbb Z_p}q^{-t}\left(\frac{1-q^{x+t}}{1-q}\right)^k d\mu_q(t)\\
&=\left(\frac{1-q^m}{1-q}\right)^k\frac{1+q}{1+q^m}\sum_{i=0}^{m-1}(-1)^i\int_{\Bbb
Z_p}\left(\frac{1-q^{m(t+\frac{x+i}{m})}}{1-q^m}\right)^kq^{-mt}
d\mu_{q^m}(t).
\endaligned\tag3$$

By (2) and (3), we see that
$$\aligned
&\left(\frac{1-q^N}{1-q}\right)^m\int_{\Bbb
Z_p}\left(\frac{1-q^{N(x+\frac{a}{N})}}{1-q^N}\right)^m q^{-Nx}
d\mu_{q^N}(x)\\
&=\frac{1+q^N}{1+q^{Np}}\sum_{i=0}^{p-1}\left(\frac{1-q^{Np}}{1-q}\right)^m(-1)^i\int_{\Bbb
Z_p}\left(
\frac{1-q^{pN(x+\frac{a+iN}{pN})}}{1-q^{pN}}\right)^mq^{-pNx}
d\mu_{q^{pN}}(x).
\endaligned\tag4$$

From (2), (3) and (4) we note that the $p$-adic  integration is
given by
$$T_q(s, a, N: q^N)= \frac{ 1+
q^N}{1+q^{pN}}\sum_{\Sb i=0\\ a+iN \not\equiv 0 (mod \, p)
\endSb}^{p-1}(-1)^i T_q(s, (a+iN)_{pN}, p^N: q^{pN})$$
such that
$$\aligned
&T_q(m, a, N: q^N)\\
&= \left( \frac{1-q^N}{1-q}\right)^m \int_{\Bbb
Z_p}\left(\frac{1-q^{N(x+\frac{a}{N})}}{1-q^N}\right)^m q^{-Nx}
d\mu_{q^N}(x)\\
&-\left(\frac{1-q^{pN}}{1-q}\right)^m\int_{\Bbb
Z_p}\left(\frac{1-q^{pN(x+\frac{(p^{-1}a)_n}{N})}}{1-q^{pN}}\right)^nq^{-pNx}d\mu_{q^{p^N}}(x),
\endaligned$$
where $(p^{-1}a)_N$ denotes the integer $x$ with  $0\leq x <N,$
$px\equiv a$ (mod $N$) and $m$ is integer with $ m+1\equiv 0$ (mod
$p-1$).

Hence, we have
$$\aligned
&\sum_{M=1}^{k-1}\left(\frac{1-q^M}{1-q}\right)(-1)^{M-1}T_q(m, hM,
k:q^k)=\left(\frac{1-q^k}{1-q}\right)^{m+1}S_{m,q}(h, k:
q^k)\\
&-\left(\frac{1-q^k}{1-q}\right)^{m+1}\left(\frac{1-q^{kp}}{1-q^k}\right)^m
S_{m,q}((p^{-1}h)_k , k: q^{pk}),
\endaligned\tag5$$
where $ p\nmid k$ and $p\nmid hM$ for each $ M$.

For $s \in \Bbb Z_p$, let us define $p$-adic Dedekind type DC sums
as follows:
$$
S_{p,q}(s: h, k: q^k)=
\sum_{M=1}^{k-1}\left(\frac{1-q^M}{1-q}\right)(-1)^{M-1}T_q(s, hM,
k: q^k).$$

 Then there exists a
continuous function $S_{p,q}(s:h, k:q^k) $ on $\Bbb Z_p$ which
satisfies
$$\aligned
&S_{p,q}(m:h, k:q^k)=\left(\frac{1-q^k}{1-q}\right)^{m+1}S_{m,q}(h,
k:q^k)\\
&-\left(\frac{1-q^k}{1-q}\right)^{m+1}\left(\frac{1-q^{kp}}{1-q^k}\right)^{m}S_{m,q}(
(p^{-1}h)_k, k: q^{pk}), \text{ where $ m+1\equiv 0$ (mod $p-1$)}.
\endaligned$$

Remark. Note that
$$ S_m (h,k: q^l)=\sum_{M=1}^{k-1} \left(\frac{1-q^M}{1-q^k}\right)(-1)^{M-1}E_{m,q^l}(\{ \frac{hM}{k}\}).$$
It is easy to see that $S_{p,1}(s: h, k:1)$ is the $p$-adic analogue
of higher order Dedekind type DC sums $k^mS_{m+1}(h, k)$.

\vskip 20pt

\quad Taekyun Kim

\quad Division of General Education-Mathematics, Kwangwoon
University, Seoul 139-701, S. Korea \quad e-mail:\text{
tkkim$\@$kw.ac.kr}

\enddocument